\documentclass[12pt]{amsart}
\usepackage{comment}

\usepackage{amssymb,amsfonts,amsmath,amsthm,stackrel,enumerate,mathtools,soul}

\usepackage[T1]{fontenc}
\usepackage{lmodern}
\usepackage[letterpaper, margin=1in]{geometry}

\usepackage{enumerate}
\usepackage{hyperref}
\usepackage{multicol}
\usepackage{mathrsfs}  
\usepackage[dvipsnames]{xcolor}
\usepackage{graphicx,float,tikz,pgfplots,tikz-cd,subcaption}
\usetikzlibrary{fit,matrix,arrows, positioning,cd,
	intersections,patterns,pgfplots.fillbetween,decorations}
\usepgfplotslibrary{fillbetween}
\pgfplotsset{compat=1.12}
\usepackage[margin = 2cm]{caption}\usepackage{scalerel}
\usepackage[usestackEOL]{stackengine}

\def\dashint{\,\ThisStyle{\ensurestackMath{%
            \stackinset{c}{.2\LMpt}{c}{.5\LMpt}{\SavedStyle-}{\SavedStyle\phantom{\int}}}%
        \setbox0=\hbox{$\SavedStyle\int\,$}\kern-\wd0}\int}

\newcommand{\Z}{\mathbb Z}

\newcommand{\R}{\mathbb R}
\newcommand{\C}{\mathbb C}

\newcommand{\calB}{\mathcal{B}}

\newcommand{\calF}{\mathcal{F}}
\newcommand{\calH}{\mathcal{H}}

\newcommand{\loc}{\text{loc}}

\renewcommand{\phi}{\varphi}

\newcommand{\eps}{\varepsilon}
\newcommand{\e}{\varepsilon}

\newcommand{\parens}[1]{\left({#1}\right)}
\newcommand{\cb}[1]{\left\{{#1}\right\}}

\newcommand{\pder}[2]{\frac{\partial {#1}}{\partial {#2}}}

\newcommand{\abs}[1]{\left|{#1}\right|}
\newcommand{\ip}[2]{\left\langle{#1},{#2}\right\rangle}

\newcommand{\vol}{\operatorname{vol}}
\newcommand{\dist}{\operatorname{dist}}
\newcommand{\cone}{\mathcal C_\alpha}
\newcommand{\hcone}{\hat{\mathcal C}_\alpha}
\newcommand{\hu}{\hat{u}}

\newtheorem{Thm}{Theorem}[section]
\newtheorem{Prop}[Thm]{Proposition}
\newtheorem{Lem}[Thm]{Lemma}

\theoremstyle{definition}
\newtheorem{Def}[Thm]{Definition}

\theoremstyle{remark}
\newtheorem{Rem}[Thm]{Remark}
\numberwithin{equation}{section}
\numberwithin{figure}{section}

\author{Christian Cofoid}
\address{Department of Mathematics, Indiana University,
Bloomington, IN 47405}
\email{christian.a.cofoid@protonmail.com}
\author{Dmitry Golovaty}
\address{Department of Mathematics, The University of Akron, Akron OH 44325}
\email{dmitry@uakron.edu}

\author{Etienne Sandier}
\address{D\`epartement de Math\'ematiques
Universit\'e Paris 12 
94010  Cr\'eteil Cedex, France
}
\email{sandier@u-pec.fr}

\author{Peter Sternberg}
\address{Department of Mathematics, Indiana University,
Bloomington, IN 47405}
\email{sternber@iu.edu}
\title{A Ginzburg-Landau problem on a circular cone}

\begin{document}

\begin{abstract}
We carry out an asymptotic analysis for a Ginzburg-Landau type model for tangent vector fields defined on a cone. The results, in the spirit of \cite{BBH}, establish the degree and asymptotic location of vortices, one of which must be situated at the tip of the cone.
\end{abstract}

\maketitle
\section{Introduction}
With an eye towards capturing the behavior of nematic liquid crystals deposited on a singular surface, we undertake an asymptotic analysis of minimizers of a Ginzburg-Landau type energy defined for tangent vector fields. There is a vast literature on mathematical modeling of nematics on curved {\it smooth} surfaces, using, for example, variants of Ginzburg-Landau, Oseen-Frank and Landau-deGennes, with some considering intrinsic energy and other extrinsic, including e.g. \cite{CanSeg,IJ21,NV12,RVK,SSV} to name but a few. Here our interest lies in understanding how the presence of a singularity on the surface will effect the morphology of minimizing configurations of nematics, and so we pursue the question for tangent vector fields defined on the surface of a cone. 

Our primary goals are two-fold: (1) to understand the effect of a singularity on the asymptotic location of vortices and resulting distribution of degrees, and (2) to explore how well the sizable arsenal of techniques developed over the years in such works as \cite{BBH,Jer,San98,SS07} can be adapted to handle the presence of these singularities. Regarding the first goal, we mention the recent articles \cite{Long_Nelson,ZN} where physicists analyze a model, not of Ginzburg-Landau type, for nematics deposited on the surface of a cone. They note the emergence of `fractional charges' induced by a singularity at the tip of a cone as well as the effect of the opening angle of the cone on defect behavior. Using a Ginzburg-Landau type model, our work will also reveal the inevitable emergence of fractional degree and a dependence of vortex behavior on cone angle.

To be precise, we let $\cone\subset \R^3$ be a finite circular cone with generator of length $1$, and
 with its tip at the origin removed. We endow $\cone$ with the metric induced by $\R^3$, and we note that $\cone$ can be unwrapped isometrically along a generator into a circular sector   $\hcone$ in the plane with sides identified. The parameter $\alpha$ denotes the opening angle $\alpha$ of this planar sector. We will denote this isometry by $\mathcal I:\cone\to\hcone$.

Denoting by $T\cone$ the tangent space to the cone, we take $u\in H^1(\cone;T\cone)$ and for $\e>0$, we define the Ginzburg-Landau energy on $\cone$:
\begin{equation}
    E_\e(u,\cone):=\int_{\cone}|Du|^2+
                \frac{1}{4\e^2}(1-\abs{u}^2)^2\,dS.\label{GLdefn}
\end{equation}
Here $D u$ denotes the covariant derivative of the vector field and $dS$ denote surface measure on the cone. As was done in \cite{ZN}, we consider only a cost measuring intrinsic deformations in our model though it would also be interesting to take an extrinsic gradient.
The central problem of this article is to obtain an asymptotic description of minimizers $\{u_\e\}$ to the problem $\inf_u E_\e(u)$, taken over the class of competitors
	\begin{align}
		H^1_g &:= \cb{u\in H^1(\cone;T\cone): u|_{\partial \cone} = g}.
	\end{align}
We assume $g\in H^{s}(\partial\cone;UT\cone)$ for some $s>1/2$,
where $UT\cone$ denotes the collection of $S^2$-valued tangent vector fields. Not surprisingly, our results will hinge on the degree $\bar{d}$ of $g$, with the notion of degree on this singular surface developed carefully below
in \eqref{eqn:surf_deg}. Somewhat more surprising, as mentioned earlier, is that our results will also depend on the opening angle $\alpha$. In our main results, expressed in Theorems \ref{leadord} and \ref{thm:1},  we establish convergence of minimizers to a canonical harmonic vector field away from a finite number of vortices whose asymptotic location is determined through minimization of a renormalized energy, much in the spirit of \cite{BBH}. One of the vortex locations is always at the tip of the cone and we determine the degree at the tip in terms of $\bar{d}$ and $\alpha$.

After describing the mathematical set-up and defining the notion of degree on a cone in Section 2, we proceed in Section 3 with an adaptation of the ball growth procedure introduced in \cite{Jer99,San98}. Here the key novelty--and difficulty--is the insistence that one of the balls is always centered at the tip of the cone, so that as it grows, it absorbs any other balls it encounters. In Section 4 we construct a test sequence to obtain an upper bound on $E_\e(u_\e,\cone)$ that is optimal up to an error of order $o(1)$. Finally, in Section 5 we combine the results of the previous two sections, along with what are now standard techniques, to obtain our main results.

\noindent{\bf Acknowlegments.} D.G. acknowledges  support by an NSF grant DMS 2106551. The research of P.S. was supported by a Simons Collaboration grant 585520 and an NSF grant DMS 2106516.

\section{Preliminaries}
In order to facilitate the calculations in the analysis to follow and to define precisely the notion of degree, we now introduce some notation and a few basic notions. Letting $r$ and $\theta$ be polar coordinates on $\hcone$ with $0\leq r\leq 1$ and $0\leq \theta\leq \alpha$, we denote by $e_r$ and $e_\theta$ an orthonormal basis on $\hcone$ induced by these coordinates. This induces an orthonormal basis $\{N,T\}$ for $\cone$,
where $N$ and $T$ are the lifting via $\mathcal{I}^{-1}$ of $e_r$ and $e_\theta$, respectively. In particular, $N$ is a unit vector pointing along the generators of the cone.

A map $u:\cone \to T\cone$ can then be expressed as \begin{equation}
    u=u^{(1)}\,T+u^{(2)}\,N,\label{conevf}
\end{equation} where $u^{(1)}$ and $u^{(2)}$ are real-valued maps. In light of the isometry $\mathcal{I}$, the image of a continuous tangential vector field $u$ under $\mathcal{I}$ induces a map $\hat{u}:\hcone\to\C\approx \R^2$ given by
\begin{equation}
 \hat{u}=u^{(1)}\,e_\theta+u^{(2)}\,e_r,\label{polarhat}   
\end{equation}
or, as a complex-valued map by
\[
\hat{u}=u^{(2)}e^{i\theta}+iu^{(1)}e^{i\theta}.
\]
Furthermore, expressing $\hat{u}$ in polar coordinates on $\hcone,$ the continuity of $u$ implies that for every $r\in (0,1]$, $\hat{u}$ satisfies the jump condition
    \begin{equation}
    \label{jumphat}
  \hat{u}(r,\alpha)=\hat{u}(r,0)e^{i\alpha}.
    \end{equation}

Again, denoting by $D$ the covariant derivative operator on $\cone$, one calculates that
	\begin{equation}
			D_N u= u^{(1)}_rN + u^{(2)}_rT,\label{dN}
	\end{equation}
and
	\begin{equation}
		D_T u=\parens{\frac{u^{(1)} + u^{(2)}_\theta}{r}}T + \parens{\frac{u^{(1)}_\theta-u^{(2)}}{r}}N,\label{dT}
	\end{equation}
    where we consider the component functions $u^{(1)}$ and $u^{(2)}$ as functions of the polar coordinates $r$
and $\theta$. In the same notation, we have
\[
\abs{Du}^2=\abs{D_Nu}^2 + \abs{D_T u}^2.
\]

Next, we follow \cite{IJ21} and define the map $i:T\cone\to T\cone$ to be a smooth bundle isomorphism such that for every $x\in \cone$, $i:T_x\cone\to T_x\cone$ is an isometry which satisfies for all $v,w\in T_x\cone$
	\begin{align}
	\label{i_isometry}
		i^2(v) = -v, && \ip{iv}{w} = -\ip{v}{iw} =\vol_{{\cone}}(v,w).
	\end{align}
Note that this bundle isomorphism behaves much like the complex number $i\in \C$.  We may relate the vector fields $N$ and $T$ by the bundle map $i$:
	\begin{align*}
		i(N) = T, && i(T) = -N.
	\end{align*}
For a vector field $u:\cone\to T\cone,$ we then can define the {\it current} of $u$ as the differential 1-form
	\begin{equation}
	\label{eqn:curr_abstr}
		j(u) := \ip{iu}{Du}.
	\end{equation}
One may express $j(u)$ in terms of the basis of 1-forms $\{dr,\,d\theta\}$ which are dual to the frame $\cb{\partial_r,\partial_\theta}$ of $T\cone$:
	\begin{equation}
		\label{eqn:curr_1_form_basis}
		j(u) = \ip{iu}{D_{\partial_r}u}dr + \ip{iu}{D_{\partial_\theta}u}d\theta,
	\end{equation}
so that for $u$ as in \eqref{conevf} one finds
	\begin{equation}
	j(u)=\parens{u^{(1)}u^{(2)}_r - u^{(2)}u^{(1)}_r}dr + \parens{u^{(1)}u^{(2)}_\theta - u^{(2)}u^{(1)}_\theta +|u|^2}\,d\theta.
		\label{eqn:j(u)incoords}
	\end{equation}
 This expression coincides with the similarly defined $1$-form $j(\hat{u})$ defined on $\hcone$, where the frame $\{dr,d\theta\}$ corresponds to standard polar coordinates on the sector.

 Turning to the relationship between the current and the modulus of the covariant derivative of a tangent vector field $u$, we note the following identity
	\begin{align}
	\label{eqn:norm_grad_decomp}
		|Du|^2 = \abs{\frac{j(u)}{|u|}}^2 + |d|u||^2,
	\end{align}
where $d$ denotes the exterior derivative, cf. \cite[(37)]{IJ21}.
In particular, we find for unit tangent vector fields $v$ that 
	\begin{align}
	\label{eqn:norm_grad_unit_vf_curr}
		|Dv|^2 = |j(v)|^2.
	\end{align}

 The 1-form $j(u)$ will also be used to define the degree of a vector field on $\cone$. To this end, let $W\subset \cone\cup\{0\}$ be a simply connected open set with $C^1$ boundary. Assuming that $u$ is non-vanishing on $\partial W$, we define 
	\begin{align}
	\label{eqn:surf_deg}
		\deg_{\cone}(u,\partial W) := \left\{\begin{array}{ll}\frac{1}{2\pi}\int_{\partial W}(\frac{j(u)}{\abs{u}^2}-j(N),&0\notin W,\\1+\frac{1}{2\pi}\int_{\partial W}(\frac{j(u)}{\abs{u}^2}-j(N)),&0\in W,\end{array}\right.
	\end{align}
In light of \eqref{eqn:j(u)incoords}, we note that
	\begin{align}
		\int_{\partial W}j(N) = 
		\int_{\partial W}d\theta = \left\{
			\begin{array}{ll}
				0, & 0\notin W,
					\\
					\alpha, & 0\in W.
			\end{array}
			\right.\label{callmeAl}
	\end{align}
Hence, it follows that $\deg_{\cone}(u,\partial W)$ is an integer, since in (\ref{eqn:surf_deg}), we eliminated the fractional winding of the standard global frame on $\cone$. We also note that an analogous formula for the degree $\deg_{\hcone}(\hat{u},\partial \hat{W})$ of the vector field $\hat{u}$ holds on $\hcone$ and so
\[\deg_{\hcone}(\hat{u},\partial \hat{W})=\deg_{\cone}(u,\partial W).\]

\begin{Rem}
    On a smooth closed surface $S$, to obtain an integer value for the degree of a unit tangent vector field $u$, one defines it through the formula
	\begin{equation}
	\label{smooth}
		\deg_{S}(u,\partial W):=\frac{1}{2\pi}\int_{\partial W}j(u) + \frac{1}{2\pi}\int_W \kappa{dS},
	\end{equation}
where $\kappa$ is the Gaussian curvature: see Section 1 of \cite{IJ21} and Section 12 of \cite{doC98}.  The surface $\cone$, however, has a singularity, and the Gaussian curvature $\kappa$ is everywhere equal to zero.  To account for the effect of the singularity in the computation of degree, we subtract $j(N)$.  For an alternative perspective, one may also arrive at \eqref{eqn:surf_deg} through a limit of a sequence of smooth approximations of the cone using \eqref{smooth}.
\end{Rem}

In light of the isometry between $\cone$ and $\hcone$, we observe that the magnitude of the gradient of the vector field $\hat{u}$ on $\hcone$ is the same as the magnitude of the covariant derivative of $u$ taken on $\cone$ and that 
\begin{equation}
    \int_{\hcone}|D\hat{u}|^2+\frac{1}{4\e^2}(1-\abs{\hat{u}}^2)^2\,dA
				= \int_{\cone}|Du|^2+
                \frac{1}{4\e^2}(1-\abs{u}^2)^2\,dS.\label{conetosector}
\end{equation}
Therefore, whenever convenient, we can study our problem on the sector $\hcone$ with the straight edges identified, rather than working on the cone.

\section{Lower Bound for $E_\eps(\cdot, \cone)$}\label{sec:LowerBound}

We now proceed to establish a lower bound for the energy,
	\begin{align}
	\label{eqn:GL_energy}
		E_\eps(u,\cone) = \frac{1}{2}\int_{\cone} |Du|^2 + \frac{1}{2\eps^2}(1-|u|^2)^2\,dS
	\end{align}
where $u\in H^1_g$.  We will adapt the vortex ball method found in \cite{SS07}. We begin with the following result which justifies the decision
 in the sequel to always place one `bad ball' with center at the origin.
We will denote by
		\[B_r(a) = \cb{x\in {\cone}:\dist_{{\cone}}(x,a) < r}\]
        an open ball of radius $r$ centered at $a\in{\cone}\cup\{0\},$ where $\dist_{{\cone}}(x,y)$ is the geodesic distance between $x,y\in\cone\cup\{0\}$. 

\begin{Lem}
\label{lem: |u|=0_at_tip}
	If $u\in H^1_g$, then $|u|$ is not bounded away from zero on any punctured ball of $\cone$ centered at the origin.
\end{Lem}
\begin{proof}
	Suppose to the contrary that for some $r>0$ there is a ball $B_r(0)\subset \cone\cup\{0\}$ such that $|u|$ is bounded away from 0 on $B_r(0)\setminus\{0\}$.  That is, suppose $\gamma:= \mathrm{essinf}_{x\in B_r(0)}|u(x)|>0$.  We will argue that $u\notin H^1_g$.  For $v := \frac{u}{|u|}$ on $B_r(0)\setminus\{0\}$, we find
		\begin{align*}
			\int_{B_r(0)}|Du|^2{dS} 
				&= \int_{B_r(0)}|u|^2|Dv|^2 + |d|u||^2{dS}
			\\
				&\ge \int_{B_r(0)} |u|^2|Dv|^2{dS}
			\\
				& = \int_{B_r(0)} |u|^2|j(v)|^2\,{dS}.
		\end{align*}
	Writing the integral over $B_r(0)$ in the coordinates $(r,\theta)$, and using the definition of degree given in \eqref{eqn:surf_deg} along with Cauchy--Schwarz and \eqref{callmeAl}, we find
		\begin{align*}
			\int_{B_r(0)}|u|^2|j(v)|^2{dS} 
				&= \int_0^r \int_{\partial B_\rho(0)}|u|^2|j(v)|^2 d\calH^1\,d\rho
			\\
				&\ge \gamma^2 \int_0^r \frac{1}{\calH^1(\partial B_\rho(p))}\parens{\int_{\partial B_\rho(0)}j(v)}^2\,d\rho
			\\
				&= \gamma^2\int_0^r \frac{1}{\calH^1(\partial B_\rho(0))} \parens{2\pi \deg_{\cone}(v,\partial B_{\rho}(0)) + \alpha-2\pi}^2\,d\rho
			\\
				&= \gamma^2\int_0^r \frac{(2\pi\deg_{\cone}(v,\partial B_{\rho}(0)) + \alpha-2\pi)^2}{\calH^1(\partial B_\rho(0))} d\rho.
		\end{align*}
	But since $\calH^1(\partial B_\rho(0)) = \alpha \rho$ for any value of $\rho\in (0,r)$, and since the value of 
		\[(2\pi\deg_{\cone}(v,\partial B_\rho(0)) +\alpha-2\pi)^2\]
	is a positive constant independent of $\rho$, we see that
		\[\frac{\gamma^2}{\alpha}\int_0^1\frac{(2\pi\deg_{\cone}(v,\partial B_{\rho}(0)) + \alpha-2\pi)^2}{\rho}\,d\rho = \infty.\]
\end{proof}

We now introduce a quantity that plays a crucial role in the asymptotic analysis of the minimal energy of $E_\e$. 

\begin{Def}\label{defmdphi}
For any $d\in \Z$ and $\alpha\in (0,2\pi)$, we define 
    \begin{equation}
    \label{newm}
    m(d,\alpha) :=\min_{\substack{d_0,d_1\in\Z\\d_0+d_1 = d}}\, \frac{2\pi}{\alpha}\left(d_0-1+\frac{\alpha}{2\pi}\right)^2 + |d_1| .\end{equation}
Then through a simple calculation one can solve this minimization problem to find that   \begin{equation} \label{mform} 
  m(d,\alpha) =   \begin{cases}
 |d| + \frac{(\alpha - 2\pi)^2}{2\pi\alpha},& \text{if $d\le 0$ and $\alpha>2\pi/3$,}\\   |d-1| + \frac{\alpha}{2\pi}, &\text{otherwise.}
\end{cases}
\end{equation}
\end{Def}\
We next note the following property of the quantity $m(d,\alpha)$.
\begin{Lem}\label{lemadditive}
    Given $\bar{d}\in\mathbb Z$ and a set of integers $\left\{d_0,\ldots,d_n\right\}\subset\mathbb Z,$ such that $\sum_{j=0}^nd_j=\bar{d},$ we have
    \begin{equation}
        \label{eq:wegave}
        m(\bar{d},\alpha)\leq m(d_0,\alpha)+\sum_{j=1}^n\left|d_j\right|
    \end{equation}
\end{Lem}
\begin{proof}
We observe that
\begin{multline*}
m(d_0,\alpha)+\sum_{j=1}^n\left|d_j\right|=\sum_{j=1}^n\left|d_j\right|+\left\{\begin{array}{ll}
    |d_0|+\frac{(\alpha - 2\pi)^2}{2\pi\alpha} , & d_0\le 0,\;\alpha>2\pi/3,\\|d_0-1|+\frac{\alpha}{2\pi}, & \mathrm{otherwise}
\end{array}\right. \\ \geq \left\{\begin{array}{ll}
    \left|\sum_{j=0}^nd_j\right|+
    \frac{(\alpha - 2\pi)^2}{2\pi\alpha},
    & d_0\le 0,\;\alpha>2\pi/3,\\\\\left|\sum_{j=0}^nd_j-1\right|+\frac{\alpha}{2\pi}, & \mathrm{otherwise}
\end{array}\right.\\  =\left\{\begin{array}{ll}
    \left|\bar{d}\right|+
    \frac{(\alpha - 2\pi)^2}{2\pi\alpha},
     \;d_0\le 0,\;\alpha>2\pi/3,\\ \left|\bar{d}-1\right|+\frac{\alpha}{2\pi},\;  \mathrm{otherwise}
\end{array}\right. =m(\bar{d},\alpha).
\end{multline*}
\end{proof}

Now we are ready to state the analogue of Theorem~3 of \cite{San98}.
\begin{Thm}
\label{thm:LB1}
Let $g:\partial\cone\to UT\cone$  be a boundary map belonging to $H^s$, for an arbitrary $s>1/2$, and  $\bar d:=\deg_{\cone}(g,\partial\cone)$. Then for all
$u:\cone\to T\cone$ such that $u=g$ on $\partial\cone$, we have
$$E_\eps(u,\cone)\ge \pi m(\bar d,\alpha) \log\frac{1}{\e} -C,$$
where $C$ depends only on $\alpha$ and $M:=\|g\|_{H^s}$, for an arbitrary $s>1/2$.

Moreover, if $E_\e(u,\cone)\le \pi m(\bar d,\alpha) \log(1/\eps) + C_1$, then there exists a constant $C'=C'(C_1,\alpha,M)$ such that
for any $\eta\in (C'\eps,1/C')$ there exists an admissible family of $K$ disjoint balls
 $\{B_0,\dots,B_K\}$ with radii less than $\eta$, with $B_0$ centered at the tip and $K\le |\bar d|+1$, 
  such that 
$$E_\eps\left(u,\cone\setminus\cup_{0=i}^K B_i\right) \le \pi m(\bar d,\alpha) \log\frac1\eta + C'.$$
\end{Thm}

The proof of the theorem closely follows the procedure in \cite{San98} that centers on the vortex ball process. This process begins by initially choosing an appropriate finite family of pairwise disjoint closed balls, and then evolving these balls to obtain a final family satisfying the conclusions of Theorem \ref{thm:LB1}.  The principal distinction with the standard ball growth procedure  is that here we use only admissible families of balls, i.e. families with a ball centered at the tip of the cone. We choose the presentation of the ball-growth process of \cite{SS07}.

\begin{Def}\label{defadmissible}
	A family of pairwise disjoint closed balls $\calB$ on $\cone\cup\{0\}$ is said to be {\em admissible} if one of its balls is centered at the origin. We will denote this ball by $B_0$.
We will also denote by $r(\calB)$ the sum of the radii of the balls comprising the collection $\calB$.
\end{Def}

 Whenever convenient, we will consider the cone $\cone$ to be extended to a larger cone without introducing any new notation. Along these lines, given any $u\in H^1_g$, when convenient we will consider $u$ to be extended outside the original cone to take the value of $g$ along any (infinite) generator of $\cone$.
\vskip.1in
We first gather some useful facts about balls on the cone ${\cone}$. A ball will be said to be {\em self-intersecting} if it intersects every ray emanating from the origin.
   Note that if $\alpha\geq\pi$ then every self-intersecting ball contains $0$. 

\begin{Lem}\label{lem:ball_stuff} We have:
\begin{enumerate}
\item Any closed self-intersecting ball $B$ in $\cone$ which does not contain $0$ is included in a closed ball $B_0$ centered at $0$ of radius $\left(1+\frac{2\pi}{\alpha}\right)r(B)$.
\item Assume $B$ and $B'$ are two intersecting closed balls on the cone ${\cone}\cup \{0\}$. Then the following holds.
	\begin{enumerate}[\quad (i) ]
		\item If neither $B$ nor $B'$ contains $0$, then there is a closed ball $B''\subset S$ such that $B\cup B'\subset B'',$ where either\\ (a) $B''$ is not self-intersecting (so that  $0\notin B''$) and  is of radius $r(B) + r(B')$  or\\ (b) $B''$ is the ball centered at $0$ of radius $\left(1+\frac{2\pi}{\alpha}\right)(r(B) + r(B'))$.
        \item If one of the balls contains $0$, then both are are included in the ball $B''$ of radius $2(r(B') + r(B))$ centered at $0$.
		\item If $B'$ is centered at $0$ then $B\cup B'$ is included in the ball $B''$ of radius $r(B') + 2r(B)$ centered at $0$.
	\end{enumerate}  
  \end{enumerate}
\end{Lem}
\begin{proof}  Item 1, which is only relevant when $\alpha<\pi$, follows by elementary trigonometry.
For Item 2, Case (i), if $B\cap B'\neq\varnothing$ then $B\cup B'$ is included in a ball $\tilde B$ of radius $r(B)+r(B')$. Then either this ball is not self-intersecting, which corresponds to case (i)(a), or it is self-intersecting. In the latter case, from Item 1 the ball $\tilde B$ is included in a ball $B''$ centered at the tip with radius $\left(1+\frac{2\pi}{\alpha}\right)(r(B) + r(B'))$. Item 2, Cases (ii) and (iii) are immediate.
\end{proof}

We use Lemma \ref{lem:ball_stuff} to establish the following analog of Theorem 4.2 of \cite{SS07}. 

\begin{Thm}[Ball Growth]\label{thm:BallGrowth}
	Assume that $\calB_0$ is an admissible collection of balls contained in $\cone\cup\{0\}$. Then there is a family $\cb{\calB(t):t\ge 0}$ of admissible collections of balls, such that $\calB(0) = \calB_0$ and 
	\begin{enumerate}
		\item For every $0\le t\le s$
		  \begin{align}\label{eqn:ball_growth_1}
                    \bigcup_{B\in \calB(t)}B\subset \bigcup_{B\in \calB(s)}B.\end{align}
		\item There is a finite set $T\subset \R_+$ such that for every $[t,s]\subset \R_+\setminus T$, we have
			\begin{align}
   \label{eqn:ball_growth_2}
                   r(\calB(s)) = e^{s-t}r(\calB(t)).
                \end{align}
		\item For every $t\in \R_+$
					\begin{align}\label{eqn:ball_growth_3}
   e^tr(\calB_0)\le r(\calB(t))\le \left(1+\frac{2\pi}{\alpha}\right)e^tr(\calB_0).\end{align}
   	\end{enumerate}
\end{Thm}

\begin{figure}[H]
\centering
	\includegraphics[width = 0.4\textwidth]{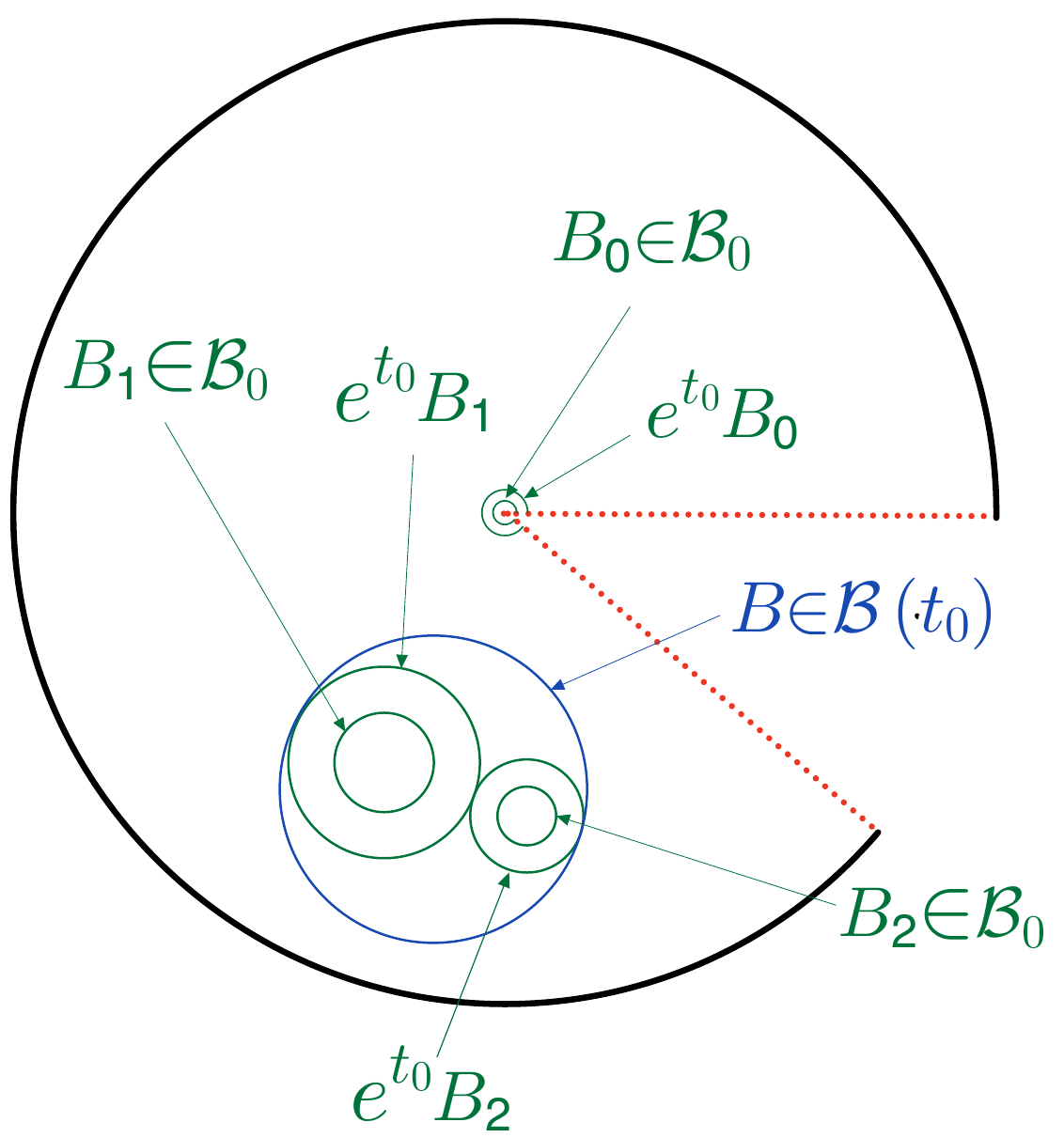}
	\caption{The merging of two balls, grown respectively from $B_0\in\mathcal{B}_0$ and $B_1\in\mathcal{B}_0$ to obtain $B\in \calB(t_0)$ of radius $r(B) = e^{t_0}\left(r(B_0) + r(B_1)\right)$. {To reduce clutter in this and subsequent figures, we drop the supercript indicating time in the notation for balls.}}
	\label{fig:merge_no_sing_pt}
\end{figure}

\begin{figure}[H]
	\centering
	\includegraphics[width = 0.4\textwidth]{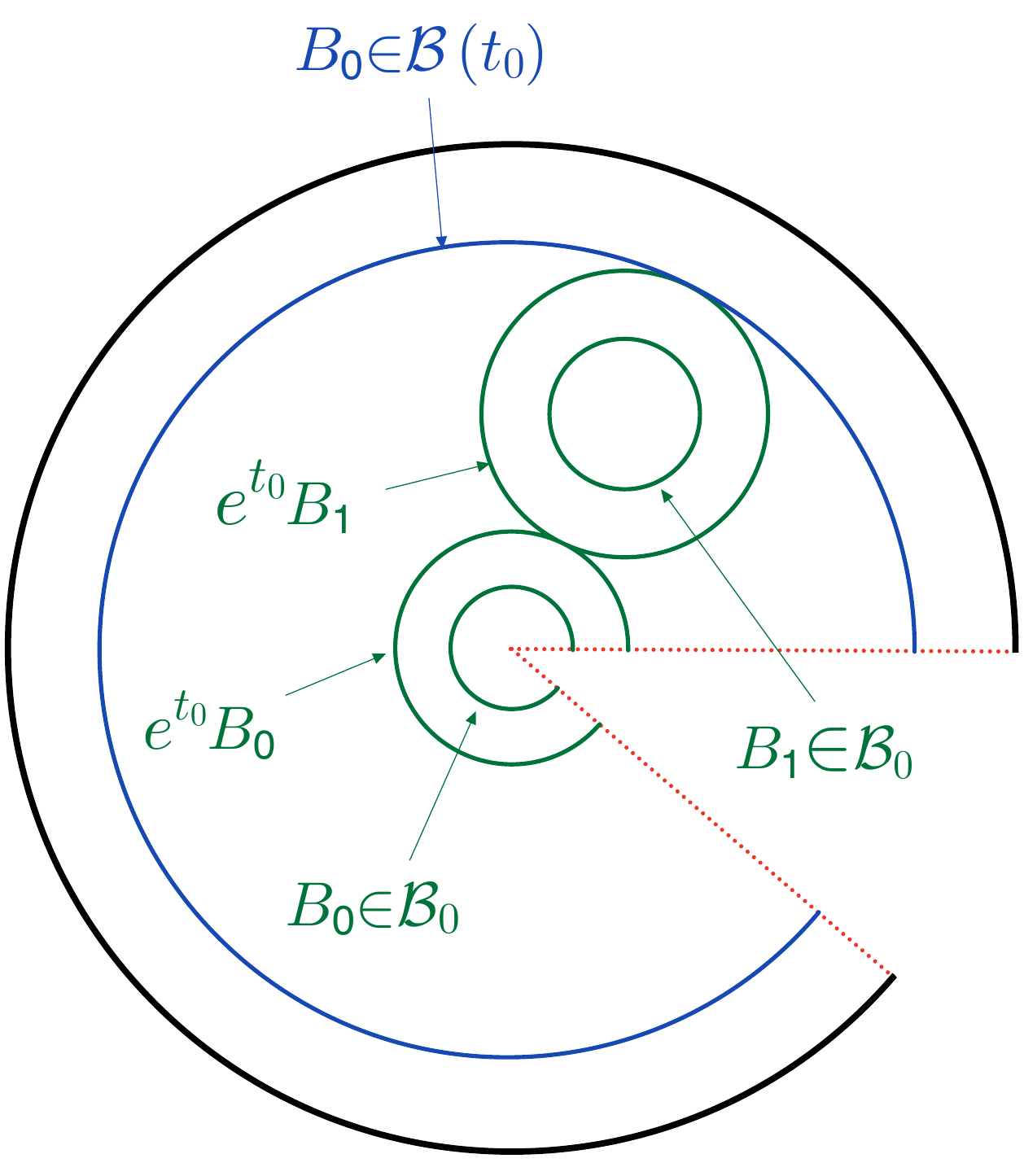}
	\caption{The merging of two balls, grown respectively from the central ball $B_0\in\mathcal{B}_0$ and $B_1\in\mathcal{B}_0$ to obtain the new central ball $B_0\in \calB(t_0)$ of radius $e^{t_0}\left(r(B_0) + 2r(B_1)\right).$ }
	\label{fig:merge_sing_pt}
\end{figure}

Figures \ref{fig:merge_no_sing_pt},  \ref{fig:merge_sing_pt}, and \ref{fig:merge_off_center_sing_pt} illustrate some of the possible scenarios for the merging process on the planar domain $\hcone.$
\begin{proof}
	We begin by growing the balls in the collection 
		\[\calB_0 = \cb{B_0,B_1,\dots, B_{k(0)}}.\]
	For $t\in[0,t_0)$, let
		\begin{align}\label{eqn:rule_of_ball_growth}
                \calB(t) = \cb{e^tB:B\in \calB_0},
            \end{align}
	 where $t_0$ is the supremum over $t\ge 0$ such that $\calB(t)$ remains pairwise disjoint and $k(t)\in\mathbb{N}$ is the number of balls in the collection $\calB(t)$ for $t>0$.  At $t_0$, at least one of the critical events described in Lemma \ref{lem:ball_stuff} occurs and we modify the collection $\calB(t_0)$ according to the procedure described in this lemma.  At $t=t_0$ this merging algorithm continues until the resulting collection of balls is disjoint. At this point the exponential ball growth resumes until another critical event is encountered. 

     \begin{figure}[H]
	\centering
	\includegraphics[width = 0.4\textwidth]{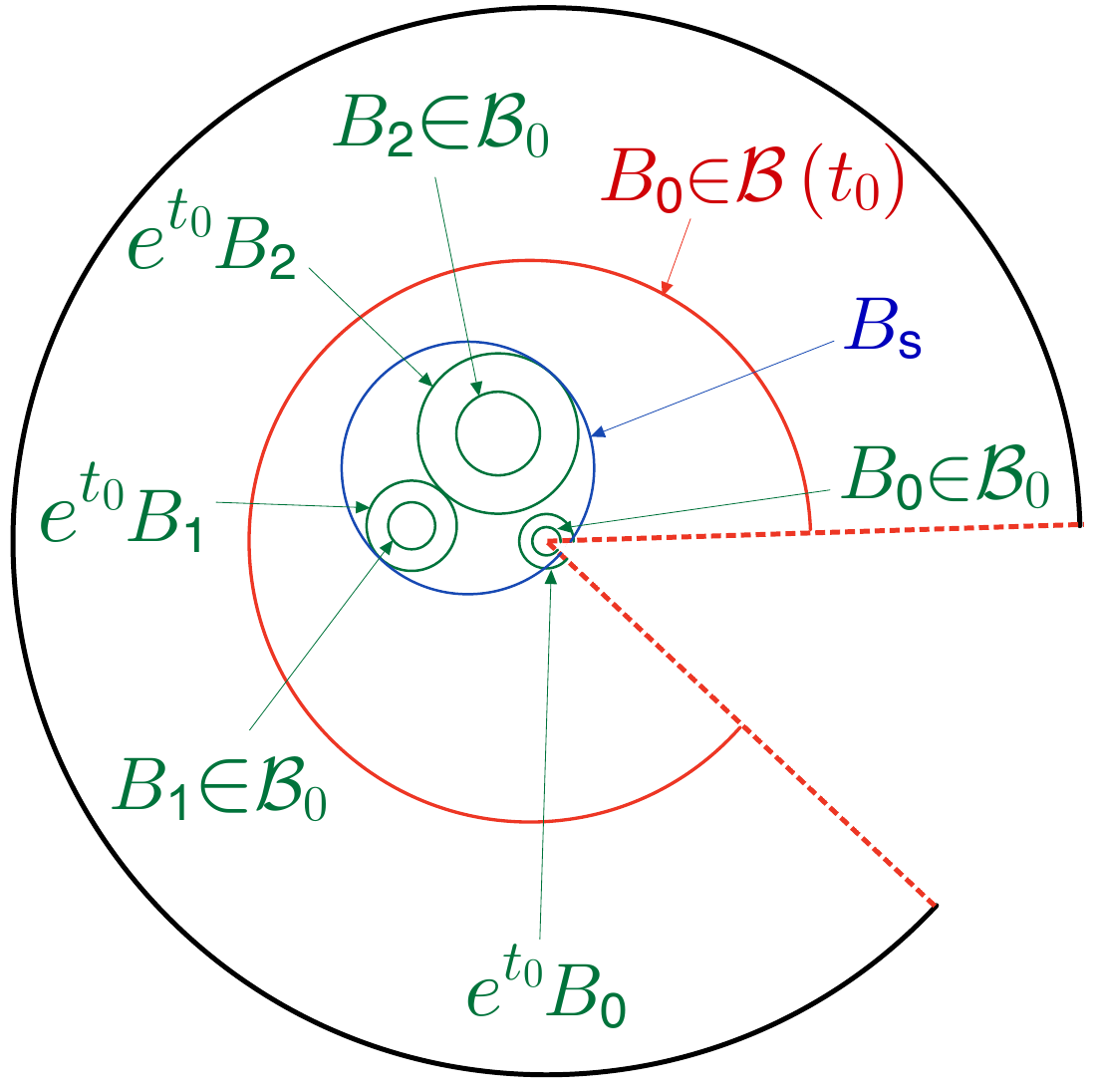}
	\caption{Growth of three balls $B_0,B_1,B_2\in\mathcal{B}_0$ generates three new balls with the larger radii and the same respective centers. The balls grown from $B_1$ and $B_2$ touch at the time $t_0$ and have to be replaced by a larger ball $B_s$ with the radius $e^{t_0}\left(r(B_1) + r(B_2)\right).$ Because $B_s$ intersects the central ball grown from $B_0$, these balls are merged in a larger new central ball of radius $e^{t_0}\left(r(B_0) + 2r(B_1) + 2r(B_2)\right).$ }
	\label{fig:merge_off_center_sing_pt}
\end{figure}

  The difference with the growth process in \cite{SS07} is that during the merging events, some of the balls see their radius multiplied by a factor which is either $2$ or $1+\frac{2\pi}{\alpha}$, hence at most $1+\frac{2\pi}{\alpha}$. As is clear from Lemma \ref{lem:ball_stuff}, this multiplication occurs to each ball at most once, when it is absorbed into the central ball. This proves \eqref{eqn:ball_growth_3}.
\end{proof}


Next, we 
recall Proposition 4.1 of \cite{SS07}. We are given a function $\calF(x,\cdot):\R_+\to \R$ which is  nondecreasing w.r.t. $r$. If $B = B_r(x)$ we write $\calF(B) = \calF(x,r)$ and if $\calB$ is a collection of balls, then $\calF(\calB)$ is the sum of $\calF(B)$ for $B\in\calB$. 
\begin{Prop}
\label{prop:calc_ineq}
	Let $\{\calB(t):t\ge 0\}$ be the collection of balls from Theorem \ref{thm:BallGrowth}, and  for each fixed $x\in {\cone}$, let $\calF(x,\cdot):\R_+\to \R$ be a nondecreasing function in $\R_+$.  Then, for $s>0$,
		\begin{align}
		\label{eqn:calc_lemma_family}
			\calF(\calB(s))-\calF(\calB_0)\ge \int_0^s \sum_{B_r(x)\in \calB(t)}r\pder{\calF}{r}\,dt.
		\end{align}
	If $B\in \calB(s)$, then
		\begin{align}
		\label{eqn:calc_lemma_single_ball}
			\calF(B)-\calF(B\cap\calB_0)\ge \int_0^s \sum_{B_r(x)\in \calB(t)\cap B}r\pder{\calF}{r}\,dt.
		\end{align}
\end{Prop}

This proposition will be used with $\calF$ being the Dirichlet energy of a unit vector field, which is bounded below in the next lemma.

\begin{Lem}(cf. Lemma 4.3 in \cite{SS07})\label{lem:amen}
	Let $B_r(x)$ be an open ball in $\cone\cup\{0\}$ such that either $0\notin B_r(x)$ or else $x=0$. Let $v:\partial B_r(x)\to UT\cone$ be a vector field of degree $d\in \Z$.  Then,
		\begin{equation}
			\frac{1}{2}\int_{\partial B_r(x)} |Dv|^2d\calH^1 \ge 
				\begin{cases}
						\dfrac{\pi}{r}d^2, & \text{if $0\notin B_r(x)$}\\
                        \dfrac{\big(2\pi(d-1)+\alpha\big)^2}{2r\alpha}, & \text{if $x=0$.}
                \end{cases}
		\end{equation}
\end{Lem}


\begin{proof}
For unit-length vector fields $v$, since $|Dv|^2=|j(v)|^2,$ cf. \eqref{eqn:norm_grad_unit_vf_curr}, using Cauchy--Schwarz we find
		\begin{align*}
			\frac{1}{2}\int_{\partial B_r(x)}|Dv|^2\,d\calH^1 
				&= \frac{1}{2}\int_{\partial B_r(x)}|j(v)|^2\,d\calH^1 
			\\
				&\ge \frac{1}{2\calH^1(\partial B_r(x))}\parens{\int_{\partial B_r(x)}j(v)}^2.
		\end{align*}
	From (\ref{eqn:surf_deg}), we also have that
		\begin{align}
		\label{eqn:int_current}
			\int_{\partial B_r(x)}j(v) = \left\{
				\begin{array}{cl}
					2\pi d & \mbox{if}\;0\notin B_r(x)
					\\
					2\pi(d-1)+\alpha & \mbox{if}\;0\in B_r(x)
				\end{array}
			\right..
		\end{align}
The conclusion of the lemma follows from
\[\calH^1(\partial B_r(x))=\left\{\begin{array}{ll}
     2\pi r,& \mbox{if}\;0\notin B_r(x),  \\
      r\alpha,& \mbox{if}\;x=0. 
\end{array}\right.\]
\end{proof}

Now, we present an analogous proposition to Theorem~1 of \cite{San98} or to Proposition 4.2 of \cite{SS07}. For this proposition it will be convenienct to denote the $j^{th}$ ball in a collection $\calB(t)$ by $B_j^{(t)}$.

\begin{Prop}
\label{prop:US_lb}
	Let $\omega$ be the compact set $\omega = \bigcup_{B\in \calB_0}B$ obtained by taking the union over an admissible initial family $\calB_0$ of balls. For $v:\cone\setminus \omega\to UT\cone$,  if $\cb{\calB(t)}_{t\in \R_+}$ is the family of admissible collections of balls grown from $\calB_0$ according to the procedure in Theorem \ref{thm:BallGrowth}, and $B\in \calB(s)$ for some $s>0$, then
 
 \begin{equation}
 \label{eqn:uvf_lb_p_or_no_p}
 E_\eps(v,B\setminus \omega)\ge
 \begin{cases}
 \pi|d_B|\parens{\log \frac{r(\calB(s))}{r(\calB_0)}-\log 2}, & \text{if $B\neq B_0^{(s)}$}\\
\pi m(d_B,\alpha)\parens{\log \frac{r(\calB(s))}{r(\calB_0)} - \log 2}, & \text{if $B= B_0^{(s)}$,}
\end{cases}
\end{equation}
	where $m(d_B,\alpha)$ is defined in Definition~\ref{defmdphi} and $d_B:={\rm{deg}}_{\cone}(v,\partial B)$.
\end{Prop}

\begin{proof}
	We define the function $\calF:\cone\times \R_+\to \R$ to be
		\[\calF(x,r) := \frac{1}{2}\int_{B_r(x)}|Dv|^2\,dS.\]
	Using Proposition \ref{prop:calc_ineq}, and keeping in mind the relationship between the radii of growing balls and time, we have that for any $B\in \calB(s)$,
		\begin{equation}\label{FminusF}
			\calF(B)-\calF(\calB_0\cap B) 
				\ge \int_0^s \sum_{{B_j^{(t)}\subset B}}r \pder{\calF}{r}\,dt= \int_0^s \sum_{{B_j^{(t)}\subset B}}\frac{r}{2}\int_{\partial B_j^{(t)}}|Dv|^2\,d\calH^1\,dt.
		\end{equation}
	If $B\neq B_0^{(s)}$, then because
		\begin{align}
			\sum_{{B_j^{(t)}\subset B}}\left|d_{B_j^{(t)}}\right|^2 \ge \sum_{{B_j^{(t)}\subset B}}\left|d_{B_j^{(t)}}\right| \ge |d_B|,
		\end{align}
	one finds through an application of Lemma \ref{lem:amen} that
		\begin{align}
			\int_0^s \sum_{{B_j^{(t)}\subset B}}\frac{r}{2}\int_{\partial B_j^{(t)}}|Dv|^2\,d\calH^1\,dt \ge \pi\int_0^s \sum_{{B_j^{(t)}\subset B}}|d_{B_j^{(t)}}|^2\,dt \ge \pi s|d_B|.
		\end{align}
	But we also know from \eqref{eqn:ball_growth_3} that 
		\begin{equation}
			\log \frac{r(\calB(s))}{r(\calB_0)} - \log 2 \le s \le \log \frac{r(\calB(s))}{r(\calB_0)}\label{lastlab}
		\end{equation}
since $\alpha\in (0,2\pi)$,	and so 
		\begin{align}
			\calF(B)-\calF(\calB_0\cap B)\ge \pi |d_B|\parens{\log \frac{r(\calB(s))}{r(\calB_0)} -\log 2}.
		\end{align}
	Next, if $B=B_0^{(s)}$, then for each $t\in (0,s)$, there a single ball $B_0^{(t)}\in \calB(t)\cap B$ which is centered at $0$.  Therefore, applying again Proposition \ref{prop:calc_ineq} and Lemma~\ref{lem:amen} to \eqref{FminusF}, it follows that
$$\calF(B)-\calF( \calB_0\cap B) \ge \int_0^s\left[\sum_{\substack{{B_j^{(t)}\subset B}\\j\neq0}}\pi|d_{B_j^{(t)}}|^2 +
\frac{\left(2\pi(d_{B_0^{(t)}} -1)+ \alpha\right)^2}{2\alpha}
\right]\,dt.$$
Denoting by $d_1(t)$ the sum of the degrees of the balls in $\calB(t)$ which are included in $B$ and do not contain $0$ we thus have 
$$\calF(B)-\calF(B\cap \calB_0) \ge \pi\int_0^s |d_1(t)| +
\frac{\left(2\pi(d_{B_0^{(t)}} -1)+ \alpha\right)^2}{2\pi\alpha}
\,dt\ge \pi s\, m(d_B,\alpha),
$$
where we have used the fact that $d_1(t)+d_{B_0^{(t)}} = d_B$ and  Definition \ref{newm}
of $m(d_B,\alpha)$. The proposition now follows from \eqref{lastlab}.
\end{proof}

\begin{proof}[Proof of Theorem \ref{thm:LB1}] In light of Proposition \ref{prop:US_lb}, the proof follows that of Theorem~3 in \cite{San98}, the only difference being that in our case, the radius $r(K)$ of a compact set $K$ is defined as the infimum of the sum of the radii of the (intrinsic) balls from an {\em admissible} collection covering $K$. We must then check that if $u$ is a, say, $C^1$ tangent vector field on $\cone$, which is of modulus $1$ on $\partial\cone$, then for any $0<t<1$ which is a regular value of $|u|$, we still have for this definition of the radius that 
\begin{equation}\label{truc}
    r(\Omega_t) \le C \text{length}(\gamma_t),
\end{equation} 
where $\Omega_t = \{|u|\le t\}$ and $\gamma_t = \{|u| = t\}$.

This follows from the fact that, since $|u| = 1$ on $\partial \cone$, then the boundary of  $\Omega_t$ is precisely $\gamma_t$, and therefore there exists a family of balls which cover $\gamma_t$ with total radius bounded by the length of $\gamma_t$. We may then include the balls containing the tip of the cone in a ball centered at the tip, and merge the intersecting or self-intersecting balls in the resulting collection. From Lemma~\ref{lem:ball_stuff}, it follows that the total radius of the collection is at most multiplied by a constant depending on $\alpha$ during this process. The resulting collection of closed disjoint balls contains $\gamma_t$, hence $\Omega_t$, which proves \eqref{truc}.

The rest of the proof being identical to that of Theorem~3 in \cite{San98}, we omit it.
\end{proof}

\section{Upper Bound for $\min_u E_\eps(u,\cone)$}
We now construct an upper bound for the infimum of $E_\eps$ over $H^1_g$.
 Our strategy will be to first construct a sequence of  complex-valued maps on the unit disc $B_1\subset\R^2$ having vortices of prescribed degree.  Away from small discs centered at the vortices with $\e$-dependent radii, this map will take the form of a standard `canonical harmonic map,' much in the spirit of \cite{BBH}, where we have chosen an optimal distribution of degrees so as to achieve the infimum given by $m(\bar{d},\alpha)$, cf. \eqref{mform}. Then we will generate a sequence of complex-valued maps on the sector $\hcone$ using the conformal map $P:B_1\to\hcone$ given by
\begin{equation}
    P(z)=z^{\alpha/2\pi},\label{cm}
\end{equation}
where we will also have to rotate the map to achieve the jump condition \eqref{jumphat}, see \eqref{fhat} below.
Finally we will utilize the inverse map of the isometry $\mathcal{I}:\cone\to\hcone$ to generate a sequence of tangent vector fields on $\cone$, yielding an upper bound on the minimal energy that, when paired with a corresponding lower bound, will prove to be optimal up to $O(1)$ terms. 

Lifting this construction to the cone $\cone$ will yield a sequence of maps which outside of slightly distorted balls on the cone agree with a canonical harmonic unit tangent vector field whose definition is given in the following:

\begin{Def}
    By a canonical harmonic unit tangent vector field $u_*$ on $\cone$  we mean one which, by unfolding the cone via $\mathcal{I}$, corresponds to a unit norm, complex-valued function  $\hat u_*$ on $\hcone$, smooth away from the origin and a finite number of points $\hat{b}_1,\dots \hat{b}_K$ on $\hcone\setminus \{0\}$,  having specified degrees $d_0,\ldots,d_K$. In a neighborhood of each $\hat{b}_j$, $j=1,\ldots,K$, the map $\hat u_*$ takes the form $e^{i(d_j\hat{\theta}_j+\hat{\phi}_j)}$ for a harmonic function $\hat{\phi}_j$ smooth in a neighborhood of $\hat{b}_j$, where $\hat{\theta}_j$ is the polar angle centered at $\hat{b}_j$. Moreover, in a neighborhood of the origin, $\hat u_*$ takes the form $e^{i[(1+2\pi (d_0-1)/\alpha)\hat{\theta}+\hat{\phi}_0]}$ for a harmonic function $\hat{\phi}_0$ smooth near the origin.
\end{Def}

\subsection{Definition of the sequence yielding the asymptotic upper bound}
Before we begin the construction, we point out the relation between the degree $\bar{d}:=\deg_{\cone}(g,\partial \cone)$ of the given Dirichlet condition $g:\partial \cone\to T\cone$ satisfying $\abs{g}=1$ and
that of a corresponding $S^1$-valued map on $\partial B_1$. Recalling the notation introduced in \eqref{conevf}-\eqref{polarhat},
if we write $g$ as
\[g=g^{(1)}\,T+g^{(2)}\,N,\]
and introduce $\hat{g}:\partial \hcone\to S^1$ via
\[
\hat{g}:=g^{(1)}e_\theta
+g^{(2)}e_r,\]
then in light of \eqref{jumphat}, the map, say $\tilde{g}:\partial B_1\to S^1$, given by
\[
\tilde{g}(e^{i\theta}):=\hat{g}(e^{i\alpha\theta/2\pi})e^{-i\alpha\theta/2\pi}
\]
will be continuous. Furthermore, 
in light of our definition \eqref{eqn:surf_deg} of degree on $\cone$ and  the identity \eqref{eqn:j(u)incoords}, one has that
\begin{equation}
 \deg(\tilde{g},\partial B_1)=\deg_{\cone}(g,\partial \cone)-1=\bar{d}-1.
 \label{theone}
\end{equation}

With these preliminaries recorded we can now proceed to define an $H^1(B_1;\C)$ map which through the use of $P$ and ${\mathcal{I}}^{-1}$ will lead to a sequence of tangent vector fields on $\cone$ that provide the desired upper bound on $E_\e$. 

To begin, we excise  $K$ balls from $B_1$ centered at $z_0,z_1,\ldots,z_K\subset B_1$ where $z_0=0$, $K$ is either $\abs{\bar{d}-1}$ or $\abs{\bar{d}}$, depending on $\bar{d}$ and $\alpha$ in a manner to be determined shortly, 
and where the radius $r_j^\e$ of the $j^{th}$ ball is given by
\begin{equation}
    r_j^\e:=\left\{\begin{matrix}
        \sqrt{\eps}/\abs{P'(z_j)}&\mbox{for}\;j=1,\ldots,k\\
        {\e}^{\pi/\alpha}&\mbox{for}\;j=0.
    \end{matrix}\right.\label{radeps}
\end{equation}
This choice of radii ensures that under the concatenation of $P$ and ${\mathcal{I}}^{-1}$, the image on $\cone$ of the $j^{th}$ ball will only differ from a ball centered at ${\mathcal{I}}^{-1}(P(z_j))$ of radius $\sqrt{\e}$ by an amount of order $\e^{3/4}$. Actually, in the case $j=0$, the image of $B(0,r_0^\e)$ under ${\mathcal{I}}^{-1}\circ P$ will be exactly a ball centered at the tip of the cone of radius $\sqrt{\e}$. 

For the construction of the sequence, say $\{v_\e\}$, on $B_1^\e:=B_1\setminus \cup_jB(z_j,r_j^\e)$ we will follow \cite{SaSo}. Recalling from \eqref{theone} that the total degree of all singularities must be $\bar{d}-1$, one finds that the optimal distribution of vortices--that is, the one that achieves the lower bound established in Theorem \ref{thm:LB1}--is given for $\bar{d}\not=1$ by a dichotomy depending on $\bar{d}$ and $\alpha$ as follows:
\vskip.1in
\noindent{\bf Case 1}: If either $\bar{d}>1$ or if $\bar{d}<1$ and $\alpha\leq 2\pi/3$, then we place singularities at $\abs{\bar{d}-1}$ points $z_1,\ldots,z_{\abs{\bar{d}-1}}$ {\it away} from the origin, each carrying degree equal to the sign of $(\bar{d}-1)$.\\
\noindent{\bf Case 2}: On the other hand, if $\bar{d}<1$ and $\alpha > 2\pi/3$, then we place  singularities at points  $z_0,\ldots,z_{\abs{\bar{d}}}$ of degree $-1$, where $z_0$ is the origin.

Finally, we have
\vskip.1in\noindent
{\bf Case 3}: If $\bar{d}=1$ so that ${\rm deg}(\tilde{g},\partial B_1)=0$, then we do not place any vortices in $B_1$. However, due to the singularity of the cone, we will still need to cut out a ball of radius $r_0^\e$ centered the origin of $B_1$.
On the complement of this ball, we define our construction as $v_\e =e^{i\phi}$ where $\phi$ solves
 \begin{equation}
     \Delta \phi=0\;\mbox{in}\;B_1,\quad\mbox{and}\quad e^{i\phi}=\tilde{g}\;\mbox{on}\;\partial B_1.\label{novortices}
 \end{equation}

For cases 1 and 2, so that $\bar{d}\not=1$, as is standard, we take $v_\e$ on $B_1^\e$ to be the $S^1$-valued `canonical harmonic map' associated with degree $1$ or degree $-1$ singularities at the $z_j's$, depending on the sign of $\bar{d}-1$, and depending on the Dirichlet condition $\tilde{g}$. Then writing $v_\e:B_1^\e\to S^1$ as $e^{i\Phi}$ or $e^{-i\Phi}$, again depending on the sign of $\bar{d}-1$, one has that $\Phi$ is harmonic with degree $1$ around each excised ball. Furthermore, if one introduces the harmonic conjugate $\Psi:B_1^\e\to\R$ of $\Phi$, then $\Psi$ can be expressed as either
\[
\Psi=\sum_{j=1}^{\abs{\bar{d}-1}} G_{\phi_{\tilde{g}}}(z,z_j)\quad\mbox{or}\quad
\Psi=\sum_{j=0}^{\abs{\bar{d}}} G_{\phi_{\tilde{g}}}(z,z_j),
\]
depending on the dichotomy mentioned in the previous paragraph.
Here we follow \cite{SaSo} and introduce $\phi_{\tilde{g}}:=\frac{1}{\abs{\bar{d}-1}}\tilde{g}\wedge \tilde{g}_\tau=\frac{1}{\abs{\bar{d}-1}}\cdot$(the tangential derivative of the phase of the map $\tilde{g}$), and for any $p\in B_1$, $G_{\phi_{\tilde{g}}}(z,p)$ is the Neumann Green's function with a singularity at $p$ and boundary condition $\phi_{\tilde{g}}$. Specifically, $G_{\phi_{\tilde{g}}}(z,p)$ is defined for all $(z,p)\in B_1\times B_1$ such that $z\not=p$, and  satisfies
\begin{align*}
    &\Delta_z G_{\phi_{\tilde{g}}}(z,p)=2\pi\delta_p\;\mbox{for}\;z,p\in B_1\\
    &\nabla_z G_{\phi_{\tilde{g}}}(z,p)\cdot\nu_z=\phi_{\tilde{g}}\;\mbox{for}\;z\in\partial B_1,\,p\in B_1,
    \end{align*}
    along with the convenient normalization
    \begin{equation}
        \int_{\partial B_1}G_{\phi_{\tilde{g}}}(z,p)\phi_{\tilde{g}}(z)\,ds_z=0\;\mbox{for}\;p\in B_1.
        \label{sasonormal}
    \end{equation}

Consequently, we have
\begin{equation}
 \Psi(z)=\sum_{j} \log\abs{z-z_j}+R(z,z_j),   \label{logbeh}
\end{equation}
where $R(z,p)$ is the regular part of $G_{\phi_{\tilde{g}}}(z,p)$, continuous on $\big(B_1\times\bar{B_1}\big)\cup \big(\bar{B_1}\times B_1\big),$ cf. \cite{SaSo}, Section 2. This completes the description of $v_\e$ on the punctured disc $B_1^\e$.

Next, we associate with $v_\e$ a complex-valued map on the punctured sector \[
\hcone^{\e}:= \left\{\hat{z}\in\hcone:\;\hat{z}\notin\cup_{j\geq 0}P\big(B(z_j,r_j^\e)\big)\right\}
\] via the formula 
\begin{equation}
 \hat{v}_\e(\hat{z}):=v_\e(\hat{z}^{2\pi/\alpha})e^{i\hat{\theta}},
 \label{fhat}
\end{equation}
 Here $\hat{\theta}\in [0,\alpha]$ denotes the polar angle of a point $\hat{z}$. Note that $\hat{v}_\e$, so defined, will necessarily satisfy the jump condition
\eqref{jumphat}. Before lifting this construction to the cone $\cone$ using $\mathcal{I}^{-1}$, we must provide a definition of $\hat{v}_\e$ in the holes of $\hcone^{\e}$. A minor complication arises from the fact that for $j\geq 1$, the
holes $P\big(B(z_j,r_j^\e)\big)$ are distorted discs rather than actual discs centered at $P(z_j)$. However, through a simple calculation using the definition  of $P$ and $r_j^\e$ given in \eqref{cm} and \eqref{radeps}, one finds that for $j\geq 1$ one has
\[
\abs{\hat{z}-P(z_j)}=\sqrt{\e}+O(\e^{3/4})\quad\mbox{for all}\;\hat{z}\in\partial P\big(B(z_j,r_j^\e)\big).
\]
Thus, by removing actual discs centered at each $P(z_j)$ having a slightly larger radius of size $\sqrt{\e}+O(\e^{3/4})$, and defining $\hat{v}_\e$ through \eqref{fhat} in this punctured version of the sector, we can easily remedy this complication. As the effect of this change on all energy calculations to follow will be $o(1)$, we will ignore this alteration and proceed as though each set $P\big(B(z_j,r_j^\e)\big)$ were actually a disc of radius $\sqrt{\e}$.

Again for $j\geq 1$, along $\partial P\big(B(z_j,r_j^\e)\big)$, the construction yields a degree $\pm 1$ vortex; hence in terms of a local phase variable $\hat{\theta}_j$ the phase must take the form $\pm\hat{\theta}_j+\chi(\hat{\theta}_j)$, where $\chi$ is a smooth function with $\chi(2\pi)=\chi(0)$ and $\hat{\theta}_j$ is chosen so that $\chi(0)=0=\chi(2\pi)$. One then interpolates the phase in a thin annulus of radius, say, $\e^2$ so that in this slighter smaller circle, the phase is simply $\pm\hat{\theta}$. For $j\geq 1$, one can fill in the disc with a standard radial Ginzburg-Landau vortex of degree $1$ or $-1$ as required; that is, in local polar coordinates centered at $P(z_j)$, the construction  $\hat{v}_\e(r,\hat{\theta}_j)$ takes the form $f(r)e^{\pm i\hat{\theta}_j}$ for an appropriately chosen modulus $f$. Such a construction minimizes the Ginzburg-Landau energy in a ball of radius $\sqrt{\e}$ for boundary conditions $e^{\pm i\hat{\theta}_j}$, with an energy given asymptotically by
\begin{equation}
  \pi\log\big(\frac{1}{\sqrt{\e}}\big)+\gamma+o(1),\label{gammadefn}  
\end{equation}
where $\gamma>0$ is defined through \eqref{gammadefn}, cf. \cite{BBH}, IX.1. (See \cite{Mironescu} for the optimality of the radial vortex.)

We handle the task of filling in the sector of radius $\sqrt{\e}$ centered at the tip of $\hcone$ a bit differently for two reasons. First, it requires an optimal construction in a sector rather than a disc. Second, in view of \eqref{fhat} and the fact that back in the punctured disc we will insist on either degree $0$ or degree $-1$ in a neighborhood of the origin depending on the dichotomy of Cases $1$ and $2$, we will need an optimal construction with respect to a phase on
$\{\hat{z}\in \hcone: \,\abs{\hat{z}}=\rho\}$ of the form either
$\hat{\theta}+\chi(\hat{\theta})$ or 
$(1-\frac{2\pi}{\alpha})\hat{\theta}+\chi(\hat{\theta})$ where now $\chi(\alpha)=\chi(0).$
Again, we can dispense with the complication due to $\chi$ through an $o(1)$ interpolation layer. We ignore this standard calculation in what follows.

While we strongly suspect that the optimal construction will be radial, that is, of the form $f_1(r)e^{i\hat{\theta}}$ or $f_2(r)e^{i(1-\frac{2\pi}{\alpha})\hat{\theta}}$ for some functions $f_1(r)$ and $f_2(r)$ with $r=\abs{\hat{z}}$, we will leave this question to a later investigation. For now, it will suffice to employ an approach analogous to that used in \cite{BBH} by considering the Ginzburg-Landau energy on a sector of opening angle $\alpha$ and radius $\eta$
\begin{equation}
\tilde{I}_{\e,\eta}(u):= \int_{\{\hat{z}\in \hcone:\,\abs{\hat{z}}<\eta\}}
\frac{1}{2}\abs{\nabla u}^2+\frac{1}{4\e^2}(\abs{u}^2-1)^2\,d\hat{z},\label{Iep} 
\end{equation}
subject to two possible boundary conditions corresponding to our Case 1 and Case 2 from earlier in this section. Specifically, for $u=u(r,\hat{\theta})$ we define
\begin{equation*}
\mu_1(\e,\eta):=\inf \tilde{I}_{\e,\eta}(u)\quad\mbox{for}\;u\in H^1\big(\{\hat{z}\in\hcone:\,\abs{\hat{z}}\leq \eta\}\big)
\end{equation*}
subject to boundary conditions
\begin{equation}
u(r,\alpha)=u(r,0)e^{i\alpha}\;\mbox{for}\; 0\leq r\leq \eta\quad\mbox{and}\quad u(\eta,\hat{\theta})=e^{i\hat{\theta}}\;\mbox{for}\; 0\leq \hat{\theta}\leq \alpha,\label{firstbc}
\end{equation}
and
\begin{equation*}
\mu_2(\e,\eta):=\inf \tilde{I}_{\e,\eta}(u)\quad\mbox{for}\;u\in H^1\big(\{\hat{z}\in\hcone:\,\abs{\hat{z}}\leq \eta\}\big)
\end{equation*}
subject to boundary conditions
\begin{equation}
u(r,\alpha)=u(r,0)e^{i\alpha}\;\mbox{for}\; 0\leq r\leq \eta\quad\mbox{and}\quad u(\eta,\hat{\theta})=e^{i(1-\frac{2\pi}{\alpha})\hat{\theta}}\;\mbox{for}\; 0\leq \hat{\theta}\leq \alpha.\label{secondbc}
\end{equation}
We also observe through a change of variables that 
\begin{equation}
    \mu_j\big(\e,\eta\big)=\mu_j\big(\frac{\e}{\eta},1\big)\quad \mbox{for}\; j=1,2.\label{scaling}
\end{equation}

In direct analogy with \cite{BBH}, Lemma III.1, Theorem V.3 and Lemma IX.1 we have the following result:
\begin{Prop}\label{gammaonetwo}
Defining the quantity
\[
\gamma_0(d,\alpha):=\left\{\begin{matrix}\lim_{\e\to 0}\mu_2(\e,1)-\frac{\alpha}{2}\log\big(1/\e\big)&
\text{if $d\le 0$ and $\alpha>2\pi/3$,}
\\
\lim_{\e\to 0}\mu_1(\e,1)-\frac{\alpha}{2}(1-\frac{2\pi}{\alpha})^2\log\big(1/\e\big)&\mbox{otherwise,}
\end{matrix}\right.
\]
both limits exist and so $\gamma_0(d,\alpha)$ is finite for all $\alpha\in (0,2\pi)$.
\end{Prop}
\begin{proof}
    Exactly as in the proof of \cite{BBH}, Lemma III.1 one first uses a simple scaling argument to establish the non-decreasing property of the map $t\mapsto \tilde{I}(t,1)+\frac{\alpha}{2}\log t$ for boundary conditions \eqref{firstbc} as well as for $t\mapsto \tilde{I}(t,1)+\frac{\alpha}{2}(1-\frac{2\pi}{\alpha})^2\log t$ for boundary conditions \eqref{secondbc}. Hence, the two limits in question either exist or approach $-\infty$.  Then coupling this observation with the lower bound provided by Theorem \ref{thm:LB1} and \eqref{mform}, the result follows.
\end{proof}

 On $\{\hat{z}\in \hcone: \,\abs{\hat{z}}\leq \sqrt{\e}\}$ we then take for $\hat{v}_\e$ the minimizer of the variational problem \eqref{Iep} with $\eta=\sqrt{\e}$, and boundary conditions either \eqref{firstbc} for Cases 1 and 3, or \eqref{secondbc} for Case 2.

Having completed the definition of the construction throughout the sector $\hcone$, we can now lift it via the isometric folding ${\mathcal{I}}^{-1}$ to produce a sequence, say $\{V_\e\}$, of continuous, $H^1$ tangent vector fields defined on $\cone$ expressed as in \eqref{conevf}.
\subsection{Upper bound: Calculation of the energy of the constructed sequence}

Our goal is now to compute $E_\e(V_\e,\cone)$ up to $O(1)$. To this end, we first observe that
\[
\int_{\cone}\abs{\nabla V_\e}^2\,dS=\int_{\hcone}\abs{\nabla\hat{v}_\e}^2\,
d\hat{z}, 
\]
since $\mathcal{I}$ is an isometry. But since $0$ is conformal,
we can then use \eqref{fhat} to assert that 
\begin{equation}
\int_{\cone}\abs{\nabla V_\e}^2\,dS=\int_{B_1}\abs{\nabla(v_\e e^{i\alpha\theta/2\pi})}^2\,dz.    \label{xzx}
\end{equation}
Focusing first on the energy in $B_1^\e$ where $v_\e=e^{i\Phi}$ or  $v_\e=e^{-i\Phi}$ depending on the sign of $\bar{d}-1$, we have that
\begin{align}
    \frac{1}{2}\int_{B_1^\e}\abs{\nabla(v_\e e^{i\alpha\theta/2\pi})}^2\,dz&=\frac{1}{2}\int_{B_1^\e}\abs{\nabla ({\rm{sgn}}(\bar{d}-1)\Phi+\alpha\theta/2\pi)}^2\,dz\nonumber\\
    &=\frac{1}{2}\int_{B_1^\e}\abs{\nabla\Psi}^2\,dz+{\rm{sgn}}(\bar{d}-1)\int_{B_1^\e}\frac{\alpha}{2\pi}\frac{1}{r^2}\frac{\partial \Phi}{\partial \theta}\,dz
    +\frac{1}{2}\int_{B_1^\e}\frac{\alpha^2}{4\pi^2}\frac{1}{r^2}\,dz\nonumber\\
    &=:I+II+III.\label{onetwothree}
\end{align}
To determine the asymptotic behavior of $I$ for $\e\ll 1$ one first integrates by parts: 
\[
I=\left\{\begin{matrix}-\frac{1}{2}\sum_{j=1}
^{\abs{\bar{d}-1}}\int_{\partial B(z_j,r_j^\e)} \Psi\partial_{\nu}\Psi\,ds&\mbox{if no vortex placed at the origin ({\bf Case 1})},\\
-\frac{1}{2}\sum_{j=0}
^{\abs{\bar{d}}}\int_{\partial B(z_j,r_j^\e)} \Psi\partial_{\nu}\Psi\,ds&\mbox{if one vortex placed at the origin ({\bf Case 2}),}
\end{matrix}\right.
\]
where $\partial_{\nu}\Psi$ denotes the outer normal derivative with respect to the ball $B(a_,,r_j^\e)$. We note that the use of the normalization \eqref{sasonormal} eliminates the boundary integral around $\partial B_1$.

Then, as is by now standard, cf.\cite{BBH}, one uses \eqref{logbeh} to capture all terms up to $O(1)$. The only wrinkle here is that the radii of the balls are given by \eqref{radeps}. Dropping terms that only contribute $o(1)$ to the energy, one finds:
\begin{equation*}
    I\sim \sum_{j}\int_{\partial B(z_j,r_j^\e)}\left\{
    \sum_{k}\bigg(\log\abs{z-z_K}+R(z,z_K)\bigg)\sum_{\ell}\left(\frac{1}{\abs{z-a_\ell}}\right)\right\}\,ds.
\end{equation*}
Proceeding as in \cite{SaSo}, pg. 175, when no vortex is placed at the origin one calculates that
\begin{eqnarray}
  & I \sim \pi\left(\abs{\bar{d}-1}\log\left(\frac{1}{\sqrt{\e}}\right)-
    \sum_{1\leq i,j\leq\abs{\bar{d}-1}, i\not=j}G_{\phi_{\tilde{g}}}(z_i,z_j)
    -\sum_{j=1}^{\abs{\bar{d}-1}}R(z_j,z_j)+\sum_{j=1}^{\abs{\bar{d}-1}}\log \abs{P'(z_j)}\right).\nonumber\\
&\label{Inotorigin}
\end{eqnarray}
However, when $\bar{d}-1<0$ and additionally one has $\alpha<\frac{2\pi}{3}$, then it turns out it will be energetically preferable to place one degree $-1$ vortex at the origin (i.e. at $z_0$) and $\abs{\bar{d}}$ degree $-1$ vortices elsewhere (i.e. at $z_1,z_2,\ldots,z_{\abs{\bar{d}}}$), in which case we find
\begin{equation}
   I \sim \pi\left(\bigg(\abs{\bar{d}}+\frac{2\pi}{\alpha}\bigg)\log\left(\frac{1}{\sqrt{\e}}\right)-
    \sum_{0\leq i,j\leq\abs{\bar{d}}, i\not=j}G_{\phi_{\tilde{g}}}(z_i,z_j)
    -\sum_{j=0}^{\abs{\bar{d}}}R(z_j,z_j)+\sum_{j=1}^{\abs{\bar{d}}}\log \abs{P'(z_j)}\big)\right).
\label{Iwithorigin2}
\end{equation}

Next we turn to the leading terms for $II$ appearing in \eqref{onetwothree}. As in the calculation of term $I$, the answer will depend upon whether or not one places a vortex at the origin.

It is convenient here to introduce $d(r):=\deg\big(e^{i\Phi},\partial B(0,r)\big)$, which is  monotone increasing, well-defined off of a set of radii of measure $O(\e)$, and non-negative, since $\Phi$ is constructed to have degree $+1$ around each $\partial B(z_j,r_j^\e)$. Then we calculate that 
\begin{eqnarray*}
II={\rm{sgn}}(\bar{d}-1)\int_{B_1^\e}\frac{\alpha}{2\pi}\frac{1}{r^2}\frac{\partial \Phi}{\partial \theta}\,dz&\sim \alpha\,{\rm{sgn}}(\bar{d}-1)\int_{r_0^\e}^1
\frac{1}{r}\left(
\frac{1}{2\pi}\int_0^{2\pi}\frac{\partial \Phi}{\partial \theta}\,d\theta\right)\,dr\\ &=\alpha\,{\rm{sgn}}(\bar{d}-1)
 \int_{r_0^\e}^1
\frac{1}{r}\,d(r)\,dr,
\end{eqnarray*}
where, in light of the bound $\abs{\nabla\Phi(x)}\leq \frac{C}{\abs{x-z_j}}$ in $B(z_j,r_j^\e)$, we may safely ignore the $o(1)$ contribution to $II$ from integration over the $O(\sqrt{\e})$ annuli that meet $\cup B(z_j,r_j^\e)$. As $d(r)$ is a step function jumping at $r$-values given by $\abs{z_j}$ up to an error of order $O(\e)$, one easily computes that if no vortex is placed at the origin, then 
\begin{equation}
    II\sim
    -\alpha\,{\rm{sgn}}(\bar{d}-1)\sum_{j=1}^{\abs{\bar{d}-1}}\log\abs{z_j},
    \label{IInotorigin}
    \end{equation}
    while if $\bar{d}-1<0$ and $\alpha<\frac{2\pi}{3}$ so that we place a degree $-1$ vortex at the origin,  then since $d(r)=-1$ on the interval $\sqrt{\e}^{2\pi/\alpha}<r<\abs{z_1}$, we find
 \begin{equation} II\sim 
- 2\pi\log\left(\frac{1}{\sqrt{\e}}\right)
+\alpha \sum_{j=1}^{\abs{\bar{d}}}\log\abs{z_j}.
 \label{IIwithorigin}
\end{equation}
Finally, we turn to term $III$ in \eqref{onetwothree}, whose value is independent of the dichotomy involving a possible degree $\pm 1$ vortex at the origin. We find
\[
III=\frac{1}{2}\int_{B_1^\e}\frac{\alpha^2}{4\pi^2}\frac{1}{r^2}\,dz=
\frac{1}{2}\int_{B_1\setminus B(0,r_0^\e)}\frac{\alpha^2}{4\pi^2}\frac{1}{r^2}\,dz-
\frac{1}{2}\int_{\cup_{j\geq 1}B(z_j,r_j^\e)}\frac{\alpha^2}{4\pi^2}\frac{1}{r^2}\,dz.
\]
Observing that the last integral is $O(\e)$, we obtain
\begin{align}
    III\sim &\frac{\alpha^2}{8\pi^2}\int_0^{2\pi}\int_{r_0^\e}^1\frac{1}{r}\,dr\,d\theta\nonumber\\
    &=-\frac{\alpha^2}{8\pi^2}\,2\pi\log r_0^\e=-\frac{\alpha^2}{4\pi}\log(\e^{2\pi/\alpha})=
    \frac{\alpha}{2}\log\left(\frac{1}{\sqrt{\e}}\right).\label{IIIvalue}
\end{align}

It remains to include the energy inside the balls $B(z_j,r_j^\e)$. For each $j\geq 1$, this is given by \eqref{gammadefn}, while for $j=0$, the cost is given up to $o(1)$ by either 
\begin{equation}
    \frac{\alpha}{2}\log\big(\frac{1}{\sqrt{\e}}\big)+\gamma_0(d,\alpha)\label{costnotip}
\end{equation}
when all vortices of degree $\pm 1$ are placed away from the tip of the cone (Case 1), or
\begin{equation}
 \frac{\alpha}{2}(1-\frac{2\pi}{\alpha})^2\log\left(\frac{1}{\sqrt{\e}}\right)+\gamma_0(d,\alpha),\label{costtip}   
\end{equation}
when we place a degree $-1$ vortex at the tip of the cone (Case 2), cf. Proposition \ref{gammaonetwo}.
 \vskip.1in
 
 Totaling up the cost of placing all vortices away from the tip (Case 1),  we collect the contributions from \eqref{gammadefn}, \eqref{Inotorigin}, 
 \eqref{IInotorigin}, \eqref{IIIvalue} and \eqref{costnotip} to find that
when either $\bar{d}-1>0$ or else $\bar{d}-1<0$ and $\alpha>\frac{2}{3}\pi$,  our construction yields
\begin{align}
E_\e\big(V_\e,\cone\big) &\sim
\pi\,m(\bar{d},\alpha)\log\left(\frac{1}{\e}\right)
+\abs{\bar{d}-1}\gamma+\gamma_0(\bar{d},\alpha)-\alpha\,{\rm{sgn}}(\bar{d}-1)\sum_{j=1}^{\abs{\bar{d}-1}}\log\abs{z_j}\nonumber\\
&-\pi 
    \sum_{1\leq i,j\leq\abs{\bar{d}-1}, i\not=j}G_{\phi_{\tilde{g}}}(z_i,z_j)
    -\pi\sum_{j=1}^{\abs{\bar{d}-1}}R(z_j,z_j)+\pi\sum_{j=1}^{\abs{\bar{d}-1}}\log \abs{P'(z_j)}+o(1).\label{firstupperbound}
\end{align}
On the other hand, when $\bar{d}-1<0$ and $\alpha<\frac{2}{3}\pi$, so that it is energetically preferable to place one degree $-1$ vortex at the origin, we invoke \eqref{gammadefn}, \eqref{Iwithorigin2}, \eqref{IIwithorigin}, \eqref{IIIvalue} and \eqref{costtip} to obtain
\begin{align}
 E_\e\big(V_\e,\cone\big) &\sim 
 \pi\,m(\bar{d},\alpha) \log\left(\frac{1}{\e}\right)+\abs{\bar{d}}\gamma+\gamma_0(\bar{d},\alpha)
+\alpha \sum_{j=1}^{\abs{\bar{d}}}\log\abs{z_j}\nonumber\\
 &-\pi\sum_{0\leq i,j\leq\abs{\bar{d}}, i\not=j}G_{\phi_{\tilde{g}}}(z_i,z_j) -\pi\sum_{j=0}^{\abs{\bar{d}}}R(z_j,z_j)+\pi\sum_{j=1}^{\abs{\bar{d}}}\log \abs{P'(z_j)}+o(1).
 \label{secondupperbound}
\end{align}
We have established:
\begin{Thm}
    The minimal energy of $E_\e$ on the cone $\cone$ subject to a Dirichlet condition $g:\partial \cone\to T \cone \cap S^2$ of degree $\bar{d}$ is bounded from above by the expression in \eqref{firstupperbound} when either $\bar{d}> 1$ or when $\bar{d}<1$ and $\alpha\geq\frac{2}{3}\pi$. If $\bar{d}<1$ and $\alpha<\frac{2}{3}\pi$, the minimal energy is bounded from above by the expression in \eqref{secondupperbound}.\\
    \vskip.05in
Finally, when $\bar{d}=1$, the minimal energy is bounded from above by
    \begin{equation}
    \frac{\alpha}{2}\log\left(\frac{1}{\e}\right)+\int_{B_1}\abs{\nabla \phi}^2\,dz+\gamma_0(1,\alpha),\label{energyzero}
    \end{equation}
    where $\phi$ is defined in \eqref{novortices}. We note that the coefficient $\frac{\alpha}{2}$ is simply $\pi\, m(1,\alpha)$, cf. \eqref{mform}
\end{Thm}

\begin{proof}
  The only part of the Theorem that remains to be checked  is \eqref{energyzero}. To this end, we return to \eqref{xzx} with $v_\e$ given by $e^{i\phi}$ from \eqref{novortices} in the punctured ball. After an application of the conformal map $P$, this leads to a definition of $\hat{v}_\e$ inside the punctured sector $\{\hat{z}\in\hcone: \,\abs{\hat{z}}>\sqrt{\e}\}$. Then $\hat{v}_\e$ is given by the minimizer of \eqref{Iep} subject to boundary conditions \eqref{firstbc} for $\eta=\sqrt{\e}$ in the region $
  \abs{\hat{z}}<\sqrt{\e}$. Again, we ignore an interpolation layer whose contribution to the energy is $o(1)$. The first term of \eqref{energyzero} arises as it does for Case 1 while term I that arose in Cases 1 and 2 is now replaced simply by the Dirichlet energy of $\phi$.
\end{proof}

We conclude this section by observing that the same procedure serves to yield the Dirichlet energy up to $o_\eta(1)$ of a canonical harmonic unit tangent vector field on $\cone\setminus \big(\{0\}\cup_{j=1}^K B(b_j,\eta)\big)$ for $\eta>0$, with singularities at the origin and points
$\mathbf b:=\{b_1,\dots,b_K\}\in\cone^{\times K}$ carrying degrees $\mathbf d:=d_0,\ldots,d_K$ satisfying $\sum_{j=0}^Kd_j=\bar{d}$ and the condition
\begin{equation}
        \label{degrees}
    m(\bar d,\alpha) = \parens{\frac{2\pi}{\alpha}\left(d_0-1+\frac{\alpha}{2\pi}\right)^2
     + \sum_{j=1}^K d_j^2};
    \end{equation}
    that is, the set of degrees  minimize the right-hand side subject to the condition $\sum_{j=0}^K d_j = \bar d$. 

To this end, we define
$W(\mathbf b;\mathbf d,\alpha)$ to be the minimum  of \begin{equation}
    \tilde{E}[v,b_0,b_1,\ldots,b_K]:=\lim_{\eta\to0}\left[\frac{1}{2}\int_{\cone\setminus\left(B_\eta(0)\cup_{j=1}^K B_\eta(b_j)\right)}|\nabla v|^2-\pi m(\bar{d} ,\alpha)\log\frac{1}{\eta}\right]\label{Wdefn}
    \end{equation}
over all canonical harmonic unit vector fields $v$  with singularities at the set of points $\{0\}\cup\mathbf b$, with corresponding (optimal) degrees $d_0,\dots,d_K$, and such that $v=g$ on $\partial\cone$. 

Then, repeating the calculation leading to \eqref{firstupperbound}--\eqref{energyzero}, but leaving out the core energy since we have excised balls about the singularities, we have
\begin{Thm}
For $j\geq 1$ we denote by $z_j$ the point on $B_1$ such that $b_j\in\cone$ is given by $b_j=\mathcal{I}^{-1}\big(P(z_j)\big)$. The quantity $W(\mathbf b;\mathbf d,\alpha)$ defined in \eqref{Wdefn} is given by: 
  \begin{align}
    W(\mathbf b;\mathbf d,\alpha)&=
-\alpha\,{\rm{sgn}}(\bar{d}-1)\sum_{j=1}^{\abs{\bar{d}-1}}\log\abs{z_j}
-\pi 
    \sum_{1\leq i,j\leq\abs{\bar{d}-1}, i\not=j}G_{\phi_{\tilde{g}}}(z_i,z_j)\nonumber\\
 &   -\pi\sum_{j=1}^{\abs{\bar{d}-1}}R(z_j,z_j)+\pi\sum_{j=1}^{\abs{\bar{d}-1}}\log \abs{P'(z_j)}\label{W1}\\&\mbox{if either}\;\bar{d}>1\;\mbox{or else}\;
    \bar{d}<1\;\mbox{and}\;\alpha>2\pi/3,\nonumber\\
 W(\mathbf b;\mathbf d,\alpha)&=   
 \alpha \sum_{j=1}^{\abs{\bar{d}}}\log\abs{z_j}
 -\pi\sum_{0\leq i,j\leq\abs{\bar{d}}, i\not=j}G_{\phi_{\tilde{g}}}(z_i,z_j)\nonumber\\ & -\pi\sum_{j=0}^{\abs{\bar{d}}}R(z_j,z_j)+\pi\sum_{j=1}^{\abs{\bar{d}}}\log \abs{P'(z_j)}\quad\mbox{if}\;\bar{d}<1\;\mbox{and}\;\alpha<2\pi/3,\label{W2}\\
 W(\mathbf b;\mathbf d,\alpha)&=\int_{B_1}\abs{\nabla\phi}^2\,dz\quad\mbox{where}\;\phi\;\mbox{is defined in}\;\eqref{novortices}\;\mbox{if}\;\bar{d}=1.\label{W3}
 \end{align}
 \end{Thm}

\section{Limiting harmonic map and renormalized energy}
\label{sec:sternberg}

\begin{Thm}
\label{leadord}
    Assume $g\in H^s(\partial\cone; UT\cone)$ for some $s>1/2$ and let $\bar d = \deg_{\cone}(g,\partial\cone)$.
Then, if $u_\eps$ denotes a minimizer of the energy $E_\eps(\cdot, \cone)$ in $ H^1_g$, 
    there exist points  $a_1,\ldots,a_K\in\cone$ with $K\le |\bar{d}|+1$
    such that $u_\eps$ subsequentially converges as $\eps\to 0$ to a unit vector field $u_*$ in $H^1_\loc\left(\cone\setminus\{a_0,a_1,\ldots,a_K\}\right)$, where $a_0 = 0$. 

    Moreover $u_*$ is a canonical harmonic unit tangent vector field, and its degrees $d_0,\dots,d_K$ about the points $\{a_0,a_1,\ldots,a_K\}$ minimize the right-hand side of \eqref{degrees} subject to the condition $\sum_{j=0}^K d_j = \bar d$. In particular, one finds that
\[
d_0=\left\{\begin{matrix}
    0&\mbox{when}\;\bar{d}\leq 0\;\mbox{and}\;\alpha>2\pi/3,\\
    1&\mbox{otherwise,}
\end{matrix}\right.
\]
with all other vortices having degree either $1$ or degree $-1$.

 \end{Thm}
\begin{proof}
The proof of the first statement follows the standard arguments to be found, e.g. in \cite{San98}.  Using Theorem~\ref{thm:LB1}, the Dirichlet energy is bounded uniformly outside arbitrarily small balls $B_0,\dots,B_K$, $K\le|\bar{d}|+1$. The $H^1_{loc}$ compactness follows, and the fact that the limit vector field has norm one is a consequence of the energy bound on  any compact subset of $\cone\setminus\{a_0,a_1,\ldots,a_K\}$.

To prove that the degrees are characterized through minimization of $m(\bar{d},\alpha)$ and that $u_*$ is a canonical harmonic vector field, one needs to bound from below the energy of $u_*$ in small annuli $B(a_k,\eta_0)\setminus B(a_k,\eta)$ around the singularities as in \cite{BBH}, with $\eta<\eta_0$. Denoting the Dirichlet energy simply by $E$, we observe from Theorem~\ref{thm:LB1} and the $H^1_{loc}$ convergence of $u_\e$ to $u^*$ that,
\begin{equation}\label{upp} E\left(u_*,\cone\setminus\cup_{j=0}^K B(a_j,\eta)\right)\le \pi m(\bar d,\alpha) \log\frac1\eta+C.\end{equation}
On the other hand, from the lower-bounds of Lemma~\ref{lem:amen}, we obtain after integration with respect to $r\in (\eta,\eta_0)$  that 
$$E\left(u_*,\cup_{j=0}^K \big(B(a_j,\eta_0)\setminus B(a_j,\eta)\big)\right)\ge 
\pi\left(\frac{2\pi}{\alpha}\left(d_0-1+\frac{\alpha}{2\pi}\right)^2 
 + \sum_{j=1}^K  {d_j}^2\right)\log\frac{\eta_0}{\eta},$$
 with $\sum_{j=0}^Kd_j=\bar{d}$.
Comparing the above inequality with \eqref{upp} and taking $\eta>0$ small enough we find the degrees must indeed minimize the right-hand side of \eqref{degrees}.

Finally, to prove that $u_*$ is a canonical harmonic vector field, we argue again as in \cite{BBH}. The fact that $u_*$ is harmonic follows for instance from the fact that it is a limit of minimizers of $E_\eps$. Then, if we consider the phase $\phi_j$ of $\hat u_*$ about the point $\hat{a}_j\in \hcone$ corresponding to $a_j$ under unfolding, we have from the degree condition that $\phi_j - d_j\theta_j$ is a well-defined harmonic function in a punctured neighborhood of $\hat{a}_j$. Moreover, from \eqref{upp}, the energy of $\phi_j$ is bounded in the punctured neighborhood. Hence the singularity is removable and $\phi_j$ is smooth. 
\end{proof}

We conclude by combining the upper and lower asymptotic bounds to obtain an energy expansion up to $o(1)$ and to determine the asymptotic location of the vortices.

\begin{Thm}
\label{thm:1} Assume $\bar{d}\not=1$.
The following expansion holds for the minimal energy: \[E_\eps(u_\eps,\cone)=\pi m(\bar{d},\alpha)\log{\frac{1}{\eps}}+K\gamma+\gamma_0(\bar{d},\alpha)+W(\mathbf a;\mathbf d)+o(1).\] $\partial\cone$,
where $W$ is defined in \eqref{Wdefn} and is given by the formulas \eqref{W1}--\eqref{W3}, with $K$ determined accordingly.

The points $\mathbf a=\{a_1,\ldots,a_K\}$ minimize $W$ over $\cone^{\times K}$ and $u_*$ minimizes $\tilde{E}[\cdot,a_1,\ldots,a_K]$ over canonical harmonic unit vector fields satisfying the boundary condition $g$ on $\partial\cone$, and singular at $0,a_1,\ldots,a_K$ with degrees $d_0,\dots,d_K$ such that $d_0+\ldots+d_K=\bar{d}$ and such that the degrees minimize the right-hand side of \eqref{degrees}.

\end{Thm}
\begin{proof}
In light of the test sequence $\{V_\e\}\subset H^1_g$ constructed in the previous section, we can assert that
\begin{equation}
    E_\eps[u_\eps]\le\pi m(\bar{d},\alpha)\log{\frac{1}{\eps}}+K\gamma+\gamma_0(\bar{d},\alpha)+W(\mathbf b;\mathbf d,\alpha)+o_\eps(1),\label{lastupper}
    \end{equation}
for any choice of $\mathbf b=(b_1,\dots,b_K)$ and optimal choice of degrees $\mathbf d$ in the sense of minimizing the right-hand side of \eqref{degrees}.

\medskip
To establish a corresponding lower bound, consider a subsequence of minimizers $\{u_\eps\}$  which converges to a canonical harmonic unit vector field $u_*$ in $C^1_{\loc}(\cone\setminus\{a_1,\dots,a_d\})$ as asserted in Theorem \ref{leadord}. In a neighborhood of $\hat a_j$ we have $\hu_* = e^{i(d_j\theta_{j}+\phi_j)}$, thus $\hu_*$ is well approximated by 
$e^{i(d_j\theta_{j}+\alpha_j)}$ on $\partial B(\hat a_j,\eta)$ if $\eta$ is small, where $\alpha_j = \phi_j(\hat a_j)$. It follows from the definition of $\gamma$ that if $j\geq 1$ one has 
$$\liminf_{\eps\to 0} E_\eps(\hat u_\eps,B(\hat a_j,\eta)) -  \pi \log\frac\eta\eps\ge \gamma  + o_\eta(1).$$

For $j=0$ we have from the definition of $\gamma_0(\bar{d},\alpha)$ that 
$$\liminf_{\eps\to 0} E_\eps(\hat u_\eps,B(0,\eta)) -  \pi \log\frac\eta\eps\ge \gamma_0(\bar{d},\alpha)  + o_\eta(1).$$
On the other hand, denoting by $E$  the Dirichlet energy, we have
\begin{align*}
 \liminf_{\eps\to 0} E_\eps(\hu_\eps,\hcone\setminus\cup_{j=0}^K B(\hat a_j,\eta)) & = 
E(\hu_*,\hcone\setminus\cup_{j=0}^K B(\hat a_j,\eta)) \\& =\pi m(\bar{d},\alpha)\log\frac1{\eta}+ \tilde{E}[u_*,a_1,\ldots,a_d] + o_\eta(1)\\
&=\pi m(\bar{d},\alpha)\log\frac1{\eta}+W(\mathbf a;\mathbf d,\alpha) + o_\eta(1).
\end{align*}
Combining these estimates we find that
$$\liminf_{\eps\to 0}\left( E_\eps(u_\eps) - \pi m(\bar{d},\alpha)\log\frac1\eps\right)\ge \gamma_0(\bar{d},\alpha) +  K\gamma + W(\mathbf a;\mathbf d,\alpha)
$$
Combining this lower bound with the upper bound \eqref{lastupper} proves the Theorem.
\end{proof}

\bibliographystyle{plain}
\bibliography{library}

\begin{thebibliography}{10}

\bibitem{BBH}
Fabrice Bethuel, Ha\"{\i}m Brezis, and Fr\'{e}d\'{e}ric H\'{e}lein.
\newblock {\em {G}inzburg-{L}andau Vortices}.
\newblock Modern Birkh\"{a}user Classics. Birkh\"{a}user/Springer, Cham, 2017.
\newblock Reprint of the 1994 edition [ MR1269538].

\bibitem{CanSeg}
Giacomo Canevari and Antonio Segatti.
\newblock Defects in nematic shells: a {$\Gamma$}-convergence
  discrete-to-continuum approach.
\newblock {\em Arch. Ration. Mech. Anal.}, 229(1):125--186, 2018.

\bibitem{doC98}
Manfredo~P Do~Carmo.
\newblock {\em Differential forms and applications}.
\newblock Springer Science \& Business Media, 1998.

\bibitem{IJ21}
Radu Ignat and Robert~L. Jerrard.
\newblock Renormalized energy between vortices in some {G}inzburg--{L}andau
  models on 2-dimensional {R}iemannian manifolds.
\newblock {\em Archive for Rational Mechanics and Analysis}, 239(3):1577--1666,
  2021.

\bibitem{Jer99}
Robert~L. Jerrard.
\newblock Lower bounds for generalized {G}inzburg--{L}andau functionals.
\newblock {\em SIAM Journal on Mathematical Analysis}, 30(4):721--746, 1999.

\bibitem{Jer}
Robert~L. Jerrard.
\newblock Lower bounds for generalized {G}inzburg-{L}andau functionals.
\newblock {\em SIAM J. Math. Anal.}, 30(4):721--746, 1999.

\bibitem{Long_Nelson}
Cheng Long and David~R. Nelson.
\newblock Liquid crystal ground states on cones with antitwist boundary
  conditions.
\newblock {\em Phys. Rev. E}, 111:025418, Feb 2025.

\bibitem{Mironescu}
Petru Mironescu.
\newblock Les minimiseurs locaux pour l'\'{e}quation de {G}inzburg-{L}andau
  sont \`a sym\'{e}trie radiale.
\newblock {\em C. R. Acad. Sci. Paris S\'{e}r. I Math.}, 323(6):593--598, 1996.

\bibitem{NV12}
Gaetano Napoli and Luigi Vergori.
\newblock Extrinsic curvature effects on nematic shells.
\newblock {\em Physical Review Letters}, 108(20):207803, 2012.

\bibitem{RVK}
Riccardo Rosso, Epifanio~G. Virga, and Samo Kralj.
\newblock Parallel transport and defects on nematic shells.
\newblock {\em Contin. Mech. Thermodyn.}, 24(4-6):643--664, 2012.

\bibitem{SaSo}
E.~Sandier and M.~Soret.
\newblock {$S^1$}-valued harmonic maps with high topological degree: asymptotic
  behavior of the singular set.
\newblock {\em Potential Anal.}, 13(2):169--184, 2000.

\bibitem{San98}
Etienne Sandier.
\newblock Lower bounds for the energy of unit vector fields and applications.
\newblock {\em Journal of Functional Analysis}, 152(2):379--403, 1998.

\bibitem{SS07}
Etienne Sandier and Sylvia Serfaty.
\newblock {\em Vortices in the Magnetic Ginzburg--{L}andau Model}, volume~70 of
  {\em Progress in Nonlinear Differential Equations and Their Applications}.
\newblock Birkh{\"a}user, 2007.

\bibitem{SSV}
Antonio Segatti, Michael Snarski, and Marco Veneroni.
\newblock Analysis of a variational model for nematic shells.
\newblock {\em Math. Models Methods Appl. Sci.}, 26(10):1865--1918, 2016.

\bibitem{ZN}
Grace~H. Zhang and David~R. Nelson.
\newblock Fractional defect charges in liquid crystals with p-fold rotational
  symmetry on cones.
\newblock {\em Phys. Rev. E}, 105(054703), 2022.

\end{thebibliography}

\end{document}